\documentclass[12pt]{article}
\usepackage[geometry]{ifsym}
\usepackage{verbatim}
\usepackage{fancyhdr}
\usepackage{graphicx}
\usepackage{mathrsfs}
\usepackage{amssymb}

\usepackage{mathrsfs}

\usepackage{pst-poly}     
\usepackage{pst-all}      

\newtheorem{pps}{Proposition}[section]
\newtheorem{cor}{Corollary}[section]
\newtheorem{lem}{Lemma}[section]
\newtheorem{thm}{Theorem}[section]
\newtheorem{cjt}{Conjecture}[section]

\newenvironment{pf}[1][Proof]{\noindent\textbf{#1.} }{\hfill\rule{1mm}{2mm}}

\makeatletter \@addtoreset{equation}{section} \makeatother

\begin{document}

\title{The Bondage Number of Mesh Networks \thanks {The work was supported by NNSF
of China (No. 11071233).}}
\author
{Fu-Tao Hu\quad Yong-Chang Cao\quad Jun-Ming
Xu\footnote{Corresponding author:
xujm@ustc.edu.cn} \\ \\
{\small Department of Mathematics}  \\
{\small University of Science and Technology of China}\\
{\small Hefei, Anhui, 230026, China} }
\date{}
\maketitle

\begin{quotation}

\textbf{Abstract}: The bondage number $b(G)$ of a nonempty graph $G$
is the smallest number of edges whose removal from $G$ results in a
graph with domination number greater than that of $G$. Denote
$P_n\times P_m$ be the Cartesian product of two paths
$P_n$ and $P_m$. This paper determines that the exact value of
$b(P_n\times P_2)$, $b(P_n\times P_3)$ and $b(P_n\times P_4)$
for $n\ge 2$.



\end{quotation}

\section{Introduction}

Throughout this paper, for the terminology and notation not defined
here, we refer the reader to \cite{x01, x03}. A graph $G=(V,E)$ is
considered as an undirected and simple graph, where $V=V(G)$ is the
vertex-set and $E=E(G)$ is the edge-set.

A nonempty subset $D\subseteq V(G)$ is said a dominating set in $G$
if every vertex in $G$ is either in $D$ or adjacent to a vertex in
$D$. The domination number $\gamma (G)$ of $G$ is the minimum
cardinality of all dominating sets in $G$. A dominating set $D$ is
said to be the minimum if $|D|=\gamma (G)$. The bondage number
$b(G)$ of a nonempty graph $G$ is the minimum number of edges whose
removal from $G$ results in a graph with larger domination number,
that is,
 $$
 b(G) =\min\{|B| : B\subseteq E(G), \gamma(G-B)>\gamma(G)\}.
 $$
A nonempty subset $B\subseteq E(G)$ is said a bondage set of $G$ if
$\gamma (G-B)>\gamma(G)$.

The concept of the bondage number is proposed by Fink et
al.~\cite{fjkr90} for an undirected graph and by Carlson and Develin
\cite{cd06} for a digraph. However the first result on bondage
numbers is obtained by Bauer et al.~\cite{bhns83}. There are many
research articles on the bondage number for undirected graphs and
digraphs (see, for example \cite{bhns83}$\sim$ \cite{chx09},
\cite{cd06} $\sim$ \cite{hx08}, \cite{ksk05,ky00}, \cite{ra08}
$\sim$ \cite{tv00}, \cite{zlm09}). In particular, Hu and
Xu~\cite{hx10} have showed that the problem determining bondage
number for general graphs is NP-hard.

Apart from its own theoretical interest, the study on the bondage
number is also motivated by the increasing importance in the design
and analysis of interconnection networks. It is well-known that the
topological structure of an interconnection network can be modeled
by a connected graph whose vertices represent sites of the network
and whose edges represent physical communication links. A minimum
dominating set in the graph corresponds to a smallest set of sites
selected in the network for some particular uses, such as placing
transmitters. Such a set may not work when some communication links
happen fault. The fault is possible in real world (hacking,
experimental error, terrorism, etc), so one needs to consider it.
What is the minimum number of faulty links which will make all
minimum dominating sets of the original network not work any more?
Such a minimum number is the bondage number, which measures the
robustness of a network with respect to link failures, wherever a
minimum dominating set is required for some applications.

Motivated by the above relevance of bondage number, one wants to
know how to compute it for a network. However, this computation is
generally difficult; no efficient algorithm has been proposed as
yet. Therefore, it is of significance to develop a technique to
determine bondage numbers for some special graphs or networks.
However, the exact value of the bondage number has been determined
for only a few classes of graphs, such as complete graphs, paths,
cycles and complete $t$-partition graphs (see, Fink {\it et
al.}~\cite{fjkr90} for the undirected cases, Huang and
Xu~\cite{hx06}, Zhang {\it et al.}~\cite{zlm09} for the directed
cases), trees (see, Bauer {\it et al.}~\cite{bhns83}, Hartnell and
Rall~\cite{hr92}, Hartnell {\it et al.}~\cite{hjvw98}, Topp and
Vestergaard~\cite{tv00}, Teschner~\cite{ts97}), de Bruijn and Kautz
digraphs (see, Huang and Xu~\cite{hx06}).

Let $P_n$ and $C_n$ be a path and a cycle of order $n$,
respectively. For the Cartesian product $G_1\times G_2$ of two
graphs $G_1$ and $G_2$, Dunbar {\it et al.}~\cite{dhtv98} determined
$b(C_n\times P_2)$ for $n\geqslant3$, Sohn, Yuan and
Jeong~\cite{syj07} determined that $b(C_n\times C_3)$ for
$n\geqslant 4$, Kang, Sohn and Kim~\cite{ksk05} determined
$b(C_n\times C_4)$ for $n\geqslant 4$, Huang ang Xu~\cite{hx08}
presented that $b(C_{5i}\times C_{5j})$ for any positive integers
$i$ and $j$; Cao, Yuan and Moo~\cite{cym10} determined $ b(C_n\times
C_5)$ for $n\geqslant 5$ and $n\not \equiv 3\,({\rm mod}\, 5)$,
but $b(C_n\times C_m)$ for $m\geqslant 6$
has been not determined as yet.

The mesh $P_n\times P_m$ is a very famous network, and its
domination number has been determined when $1\leqslant m \leqslant
6$ for many years~\cite{cch93, cch94, jk83, ks95}. However, its
bondage number has been not determined as yet. For $n=1$,
$P_1\times P_m$ is isomorphic to $P_m$, and $b(P_m)$ has been determined.
In this paper, we present the exact value of
$b(P_n\times P_2)$, $b(P_n\times P_3)$ and $b(P_n\times P_4)$
for $n\ge 2$.

The rest of the paper is organized as follows. Section 2 presents
some useful results, Section 3 determines $b(P_n\times P_2)$,
Section 4 determines $b(P_n\times P_3)$, and Section 5 determines
$b(P_n\times P_4)$. Some remarks are in Section 6, in which we
propose a conjecture: $b(P_n\times P_m)\leqslant 2$ for $m\geqslant
5$.

\section{Preliminary results}

Let $G_1=(V_1,E_1)$ and $G_2=(V_2,E_2)$ be two undirected graphs.
The {\it Cartesian product} of $G_1$ and $G_2$ is an undirected
graph, denoted by $G_1\times  G_2$, where $V(G_1\times
G_2)=V_1\times  V_2$, two distinct vertices $x_1x_2$ and $y_1y_2$,
where $x_1,y_1\in V(G_1)$ and $x_2, y_2\in V(G_2)$, are linked by an
edge in $G_1\times G_2$ if and only if either $x_1=y_1$ and
$x_2y_2\in E(G_2)$, such an edge is called a {\it vertical edge}, or
$x_2=y_2$ and $x_1y_1\in E(G_1)$, such an edge is called a {\it
horizontal edge}.
It is clear, as a graphic operation, that the Cartesian
product satisfies commutative associative law if identify isomorphic
graphs.

Throughout this paper, the notation $P_n$ denotes a path with
vertex-set $\{1,2,\cdots,n\}$. The $(n,m)$-mesh network, denoted by
$G_{n,m}$, is defined as the Cartesian product $P_n \times P_m$,
with the vertex-set $\{u_{i,j}|\ 1\leqslant i \leqslant n,
1\leqslant j \leq m\}$.

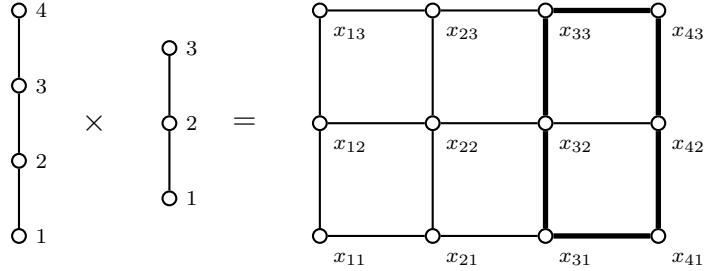
\begin{figure}[ht]
\begin{center}
\begin{pspicture}(-5,-.5)(5,3.5)

\cnode(-4,0){3pt}{1}\rput(-3.7,0){\scriptsize 1}
\cnode(-4,1){3pt}{2}\rput(-3.7,1){\scriptsize 2}
\cnode(-4,2){3pt}{3}\rput(-3.7,2){\scriptsize 3}
\cnode(-4,3){3pt}{4}\rput(-3.7,3){\scriptsize 4} \ncline{1}{2}
\ncline{2}{3} \ncline{3}{4}

\cnode(-2,.5){3pt}{1'}\rput(-1.7,0.5){\scriptsize 1}
\cnode(-2,1.5){3pt}{2'}\rput(-1.7,1.5){\scriptsize 2}
\cnode(-2,2.5){3pt}{3'}\rput(-1.7,2.5){\scriptsize 3}
\ncline{1'}{2'} \ncline{2'}{3'}

\rput(-3,1.5){$\times$}  \rput(-1,1.5){=}

\cnode(0,0){3pt}{11}\rput(.4,-.3){\scriptsize $x_{11}$}
\cnode(0,1.5){3pt}{12}\rput(.4,1.2){\scriptsize $x_{12}$}
\cnode(0,3){3pt}{13}\rput(.4,2.7){\scriptsize $x_{13}$}
\cnode(1.5,0){3pt}{21}\rput(1.9,-.3){\scriptsize $x_{21}$}
\cnode(1.5,1.5){3pt}{22}\rput(1.9,1.2){\scriptsize $x_{22}$}
\cnode(1.5,3){3pt}{23}\rput(1.9,2.7){\scriptsize $x_{23}$}
\cnode(3,0){3pt}{31}\rput(3.4,-.3){\scriptsize $x_{31}$}
\cnode(3,1.5){3pt}{32}\rput(3.4,1.2){\scriptsize $x_{32}$}
\cnode(3,3){3pt}{33}\rput(3.4,2.7){\scriptsize $x_{33}$}
\cnode(4.5,0){3pt}{41}\rput(4.9,-.3){\scriptsize $x_{41}$}
\cnode(4.5,1.5){3pt}{42}\rput(4.9,1.2){\scriptsize $x_{42}$}
\cnode(4.5,3){3pt}{43}\rput(4.9,2.7){\scriptsize $x_{43}$}

\ncline{11}{12} \ncline{12}{13} \ncline{10}{20} \ncline{11}{21}
\ncline{12}{22} \ncline{13}{23} \ncline{21}{22} \ncline{22}{23}
\ncline{20}{30} \ncline{21}{31} \ncline{22}{32} \ncline{23}{33}
\ncline[linewidth=2pt]{31}{32} \ncline[linewidth=2pt]{32}{33}
\ncline{30}{40} \ncline[linewidth=2pt]{31}{41} \ncline{32}{42}
\ncline[linewidth=2pt]{33}{43} \ncline[linewidth=2pt]{41}{42}
\ncline[linewidth=2pt]{42}{43} \ncline{40}{50} \ncline{41}{51}
\ncline{42}{52} \ncline{43}{53}

\end{pspicture}
\caption{\label{f1}\footnotesize  A $(4,3)$-mesh network
$G_{4,3}=P_4\times P_3$}
\end{center}
\end{figure}

The graph shown in Figure~\ref{f1} is a $(4,3)$-mesh network
$G_{4,3}$. It is clear, as a graphic operation, that the cartesian
product satisfies commutative associative law if we identify
isomorphic graphs, that is, $G_{n,m}\cong G_{m,n}$.

The following notations will be used in this paper. For a
positive integer $t$ with $t<n$, $G_{t,m}$ is a subgraph of $G_{n,m}$.
Denote $H_{n-t,m}=G_{n,m}-G_{t,m}$, that
is, $H_{n-t,m}$ is a subgraph of $G_{n,m}$ induced by the set of
vertices $\{u_{ij}|\ t+1\leq i \leq n, 1\leq j \leq m\}$. Clearly,
$H_{n-t,m}\cong G_{n-t,m}$. The graph shown in
Figure~\ref{f1} by heavy lines is a subgraph $H_{2,3}$ of $G_{4,3}$,
where $n=4, t=2$ and $m=3$ is a such example.

Note that both $G_{0,m}$ and $H_{n-n,m}$ are nominal graphs. For
convenience of statements, we allow $G_{0,m}$ and $H_{n-n,m}$ to
appear in this paper. If so, we specify consider that their total
dominating sets are empty.

In Addition, let $Y_i=\{u_{i,j}:j=1,2,\cdots,m\}$ for each
$i=1,2,\cdots,n$, called a set of vertical vertices of $i$ in
$G_{n,m}$.



We state some useful results on $\gamma(G_{n,m})$ to be used in this
paper.

\begin{lem}\label{lem2.1} {\rm ~\cite{jk83,ks95}}
Let $P_n$ and $C_m$ be a path and a cycle of order $n\geqslant 1$
and $m\geqslant 3$, respectively. Then

$\gamma(G_{n,2})=\lceil \frac{n+1}{2}\rceil$;

$\gamma(G_{n,3})=n-\lfloor\frac{n-1}{4}\rfloor$;

$\gamma(G_{n,4})=\left\{ \begin{array}{l}
 n+1,\ {\rm if}\ n=1,2,3,5,6\ {\rm or}\ 9  \\
n, \ ~~~~~{\rm otherwise};
\end{array}
 \right. $

$\gamma(C_m \times C_3)=n-\lfloor\frac{n}{4}\rfloor$.
\end{lem}

\begin{lem}\label{lem2.2}
Let $D$ be a dominating set of $G_{n,m}$. Then
$\gamma(G_{i,m})\leqslant |D\cap V(G_{i+1,m})|$ and
$\gamma(H_{n-i,m})\leqslant |D\cap V(H_{n-i+1,m})|$ for each
$i=1,2,\ldots,n-1$ and $m\geqslant 2$.
\end{lem}

\begin{pf}
We only need to prove that $\gamma(G_{i,m})\leqslant |D\cap
V(G_{i+1,m})|$ since $H_{n-i,m}\cong G_{n-i,m}$. Let $D'=D\cap
V(G_{i+1,m})$.

If $D'\cap Y_{i+1}=\emptyset$, then $D'$ is a dominating set of
$G_{i,m}$, and hence $\gamma_{t}(G_{i,m})\leqslant |D'|$.

Assume $D'\cap Y_{i+1}\ne\emptyset$ below. Let $B_i=\{j|\
u_{i+1,j}\in D'\}$. Then $D''=(D'\setminus Y_{i+1})\cup \{u_{i,j}|\
j\in B_i\}$ is a dominating set of $G_{i,m}$ and $|D''|\leqslant
|D'|$. Thus, we have $\gamma(G_{i,m})\leqslant |D''|\leqslant |D'|$.
The lemma follows.
\end{pf}


\section{The bondage number of $G_{n,2}$}

\begin{thm}\label{thm3.1}
$b(G_{2,2})= 3$, $b(G_{3,2})= 2$, and $b(G_{n,2})= 1$ if $n$ is odd
 and $b(G_{n,2})= 2$ if $n$ is even for $n\geqslant 4$.
\end{thm}

\begin{pf}
It is easy to verify that $b(G_{2,2})= 3$ and $b(G_{3,2})= 2$.
In the following, consider $n\ge 4$.
When $n$ is odd, we consider the domination number of
$G=G_{n,2}-u_{1,1}u_{1,2}$. Let $D$ be a minimum dominating set of
$G$. Then either $u_{1,1}\in D$ or $u_{1,2}\in D$, and either
$u_{1,2}\in D$ or $u_{2,2}\in D$. Without loss of generality we
assume that $u_{1,1},u_{1,2}\notin D$ and $u_{2,1},u_{2,2}\in D$. By
Lemma~\ref{lem2.2}, $|D\cap V(H_{n-2,m})|\geqslant
\gamma(H_{n-3,2})$. Then by Lemma~\ref{lem2.1},
 $
 |D|\geqslant 2+\gamma(H_{n-3,2}) =2+\left\lceil\frac{n-3+1}{2}\right\rceil=1+\gamma(G_{n,2}),
 $
which yields $b(G_{n,2})= 1$.

When $n$ is even, we claim that $\gamma(G_{n,2})=\gamma(G_{n,2}-e)$
for any $e\in E(G_{n,2})$.

To prove this claim, we first consider that $e$ is a vertical edge,
and let $e=u_{j,1}u_{j,2}$.

If $j$ is even, then all the vertices $u_{i,1}$, $i\equiv 1
\,({\rm mod}\,4)$, $u_{i,2}, i\equiv 3 \,({\rm mod}\,4)$, $u_{n,1}$
if $n\equiv 0\,({\rm mod}\,4)$ or $u_{n,2}$ if $n\equiv 2\,({\rm
mod}\,4)$, form a dominating set of $G_{n,2}-e$ with cardinality
$\lceil \frac{n+1}{2}\rceil$.

If $j$ is odd, then all the vertices $u_{i,1}, i\equiv 2 \,({\rm
mod}\,4)$, $u_{i,2}, i\equiv 0 \,({\rm mod}\,4)$ and $u_{2,2}$ form
a dominating set of $G_{n,2}-e$ with cardinality $\lceil
\frac{n+1}{2}\rceil$.

Assume now that $e$ is a horizontal edge. Without loss of
generality, let $e=u_{j,1}u_{j+1,1}$.

If $j\equiv 2~{\rm or}~3 \,({\rm mod}\,4)$, then all the vertices
$u_{i,1}, i\equiv 1 \,({\rm mod}\,4)$, $u_{i,2}, i\equiv 3 \,({\rm
mod}\,4)$, and $u_{n,1}$ form a dominating set of $G_{n,2}-e$ with
cardinality $\lceil \frac{n+1}{2}\rceil$.

If $j\equiv 0~{\rm or}~1 \,({\rm mod}\,4)$, then all the vertices
$u_{i,2}, i\equiv 1 \,({\rm mod}\,4)$, $u_{i,1}, i\equiv 3 \,({\rm
mod}\,4)$, and $u_{n,1}$ form a dominating set of $G_{n,2}-e$ with
cardinality $\lceil \frac{n+1}{2}\rceil$.

So we have that $b(G_{n,2})\geqslant 2$. Next, we show that
$b(G_{n,2})\leqslant 2$. Let $e_1=u_{2,1}u_{3,1}$,
$e_2=u_{2,2}u_{3,2}$, and $G'=G_{n,2}-\{e_1,e_2\}$. Then $G'$
consists of two connected components, one is $G_{2,2}$ and the other
one is $H_{n-2,2}$. By Lemma~\ref{lem2.1}, we have
 $
 \gamma(G')=\gamma(G_{2,2})+\gamma(H_{n-2,2})=2+\left\lceil
 \frac{n-2+1}{2}\right\rceil=1+\gamma(G_{n,2}),
 $
which implies $b(G_{n,2})\leqslant 2$. Thus $b(G_{n,2})=2$.
\end{pf}


\section{The bondage number of $G_{n,3}$}

\begin{pps}\label{pps4.1} {\rm (~\cite {jk83})}\
A minimum dominating set $D$ of $G_{n,3}$ is constructed as follows.
 $$
 D=\left\{\begin{array}{ll}
 \{u_{i,2}: i\equiv 1 \,({\rm mod}\,4)\}\cup \{u_{i,1}, u_{i,3}:
 i\equiv 3 \,({\rm mod}\,4)\}\ &\ {\rm if}\ n \ {\rm is\ odd},\\
 \{u_{i,2}: i\equiv 1 \,({\rm mod}\,4)\}\cup \{u_{i,1},u_{i,3}: i\equiv 3 \,({\rm
 mod}\,4)\}\cup \{u_{n,2}\}\ &\ {\rm if}\ n\ {\rm is\ even}.
 \end{array}\right.
 $$
\end{pps}

\begin{lem}\label{lem4.1}
$\gamma(G_{n,3}-u_{1,j})\geqslant
\gamma(G_{n,3})=n-\lfloor\frac{n-1}{4}\rfloor$ for each $j=1,2,3$
and $n\equiv 1,2~or~3\,({\rm mod}\,4)$.
\end{lem}

\begin{pf}
It is easy to verify that the conclusion is true for $n=1,2,3$. In
the following, assume $n\geqslant 4$. Let $G=G_{n,3}-u_{1,j}$ and
$D$ be a minimum dominating set of $G$. we only need to show
$|D|\geq n-\lfloor\frac{n-1}{4}\rfloor$.

If $(Y_1-u_{1,j})\cap D\neq \emptyset$, then $D$ is a dominating set
of $C_n\times C_3$. By Lemma \ref{lem2.1}, $
 |D|\geq \gamma(C_n\times C_3)=n-\left\lfloor \frac{n}{4}\right\rfloor
  =n-\left\lfloor\frac{n-1}{4}\right\rfloor
 $.

If $(Y_1-u_{1,j})\cap D=\emptyset$, then $|Y_2\cap D|\geqslant 2$.
By Lemma~\ref{lem2.2}, $|D\cap V(H_{n-2,3})|\geqslant
\gamma(H_{n-3,3})$. By Lemma \ref{lem2.1},
 $
 |D|\ge
 2+\gamma(H_{n-3,3})=2+n-3-\left\lfloor \frac{n-3-1}{4}\right\rfloor
 =n-\left\lfloor\frac{n-1}{4}\right\rfloor,
 $
as required.
\end{pf}

\begin{lem}\label{lem4.2}
$\gamma(G_{n,3}-u_{1,1})\geqslant
\gamma(G_{n,3})=n-\lfloor\frac{n-1}{4}\rfloor$ for $n\equiv 0\,({\rm
mod}\,4)$.
\end{lem}

\begin{pf}
Let $D$ be a minimum dominating set of $G_{n,3}-u_{1,1}$. we only
need to prove $|D|\geqslant n-\lfloor\frac{n-1}{4}\rfloor$. It is
easy to verify that the assertion is true for $n=4$. In the
following, we consider the case $n\geqslant 8$. We consider the
following three cases, respectively.

\begin{description}

\item [Case 1] $u_{1,2}\in D$ or $u_{2,1}\in D$.

In this case, $D$ is also a dominating set of $G_{n,3}$, and so
$|D|\geq \gamma(G_{n,3})=n-\lfloor \frac{n-1}{4}\rfloor$.

\item [Case 2] $u_{1,2}\notin D$, $u_{2,1}\notin D$ and $u_{1,3}\in D$.

In this case, $D\setminus \{u_{1,3}\}$ is a dominating set of
$H_{n-1,3}$ or $H_{n-1,3}-u_{2,3}$. By Lemma \ref{lem4.1},
$|D\setminus \{u_{1,3}\}|\geqslant n-1-\lfloor \frac{n-1-1}{4}
\rfloor$, and so $|D|\geqslant n-\lfloor \frac{n-1}{4}\rfloor$.

\item [Case 3] $u_{1,2}\notin D$, $u_{2,1}\notin D$ and $u_{1,3}\notin D$.

In this case, $u_{2,2},u_{2,3}\in D$. We prove the conclusion by two
subcases.

\begin{description}

\item [Subcase 3.1] $Y_3\cap D\neq \emptyset$.

Then $D\setminus \{u_{2,2},u_{2,3}\}$ is a dominating set of
$H_{n-2,3}$ or $H_{n-2,3}-u_{3,1}$ or $H_{n-2,3}-u_{3,3}$. By Lemma
\ref{lem4.1}, $|D\setminus \{u_{2,2},u_{2,3}\}|\geqslant n-2-\lfloor
\frac{n-2-1}{4}\rfloor$. Thus, $|D|\geqslant n-\lfloor
\frac{n-1}{4}\rfloor$.

\item [Subcase 3.2] $Y_3\cap D=\emptyset$.

Then $u_{4,1}\in D$.

If $u_{4,2}\in D$ or $u_{4,3}\in D$, then $D\setminus
\{u_{2,2},u_{2,3}\}$ is a dominating set of $H_{n-2,3}$ or
$H_{n-2,3}-u_{3,2}$ or $H_{n-2,3}-u_{3,3}$. By Lemma~\ref{lem4.1},
$|D\setminus \{u_{2,2},u_{2,3}\}|\geqslant n-2-\lfloor
\frac{n-2-1}{4}\rfloor$. Thus, $|D|\geqslant n-\lfloor
\frac{n-1}{4}\rfloor$.

Next, assume $u_{4,2},u_{4,3}\notin D$. Then $u_{5,3}\in D$. If
$u_{5,1}\in D$ or $u_{5,2}\in D$, then $D\setminus
\{u_{2,2},u_{2,3}, u_{4,1}\}$ is a dominating set of $H_{n-4,3}$ and
hence $|D\setminus \{u_{2,2},u_{2,3}, u_{4,1}\}|\geqslant
n-4-\lfloor \frac{n-4-1}{4}\rfloor$. Thus, $|D|\geqslant n-\lfloor
\frac{n-1}{4}\rfloor$.

If $u_{5,1},u_{5,2}\notin D$, then $D\setminus \{u_{2,2},u_{2,3},
u_{4,1},u_{5,3}\}$ is a dominating set of $H_{n-5,3}$ or
$H_{n-5,3}-u_{6,3}$. By Lemma~\ref{lem4.1}, \\
$|D\setminus \{u_{2,2},u_{2,3}, u_{4,1},u_{5,3}\}|\geqslant
n-5-\lfloor \frac{n-5-1}{4}\rfloor$. Thus, $|D|\geqslant n-\lfloor
\frac{n-1}{4}\rfloor$.
\end{description}
\end{description}
The lemma follows.
\end{pf}

\begin{cor}\label{cor3.1}
$b(G_{n,3})\leqslant 2$.
\end{cor}
\begin{pf}
By Lemma~\ref{lem4.1} and Lemma~\ref{lem4.2}, we have that
$\gamma(G_{n,3}-\{u_{1,1}u_{2,1},
u_{1,1}u_{1,2}\})>\gamma(G_{n,3})$.
\end{pf}


\begin{lem}\label{lem4.3}
$b(G_{n,3})=1$ for $n\equiv 1\ or\ 2\,({\rm mod}\,4)$ and
$n\geqslant 4$.
\end{lem}

\begin{pf}
Let $D$ be a minimum dominating set of $G_{n,3}-u_{3,1}u_{4,1}$. We
only need to prove that $|D|\geqslant 1+\gamma(G_{n,3})$ by
considering three cases, repectively.

\begin{description}

\item [Case 1] $u_{3,2}\in D$ and $u_{3,3}\in D$.

In this case, $|V(G_{3,3})\cap D|=4$. By Lemma~\ref{lem2.2}, $|D\cap
V(H_{n-3,3})|\geqslant \gamma(H_{n-4,3})$. By Lemma~\ref{lem2.1},
$|D|\geqslant 4+\gamma(H_{n-4,3})=
4+n-4-\lfloor\frac{n-4-1}{4}\rfloor=1+\gamma(G_{n,3})$.

\item [Case 2] Either $u_{3,2}\in D$ or $u_{3,3}\in D$.

In this case, $|V(G_{3,3})\cap D|=3$. Then $D'=D\setminus
V(G_{3,3})$ is a dominating set of $H_{n-3,3}$ or
$H_{n-3,3}-u_{4,2}$ or $H_{n-3,3}-u_{4,3}$. By Lemma \ref{lem4.1},
$|D|=3+|D'|\geqslant 3+ n-3-\lfloor
\frac{n-3-1}{4}\rfloor=n+1-\lfloor
\frac{n-1}{4}\rfloor=1+\gamma(G_{n,3})$.

\item [Case 3] $u_{3,2}\notin D$ and $u_{3,3}\notin D$.

In this case, $|V(G_{3,3})\cap D|=2$ or $|V(G_{3,3})\cap D|=3$. If
$|V(G_{3,3})\cap D|=3$, then $D\setminus V(G_{3,3})$ is a dominating
set of $H_{n-3,3}$. By Lemma~\ref{lem2.1}, $|D|\geqslant
3+\gamma(H_{n-3,3})=3+n-3-\lfloor\frac{n-3-1}{4}\rfloor
=1+\gamma(G_{n,3})$.

If $|V(G_{3,3})\cap D|=2$, then $V(G_{3,3})\cap
D=\{u_{1,3},u_{2,1}\}$ and $D\setminus V(G_{3,3})$ is a dominating
set of $H_{n-2,3}-u_{3,1}$. By Lemma \ref{lem4.1} or Lemma
\ref{lem4.2}, $|D|\geqslant 2+ n-2-\lfloor
\frac{n-2-1}{4}\rfloor=n+1-\lfloor
\frac{n-1}{4}\rfloor=1+\gamma(G_{n,3})$.

\end{description}
The lemma follows.
\end{pf}


\begin{lem}\label{lem4.4}
$b(G_{n,3})\geqslant 2$ for $n\equiv 0\,({\rm mod}\,4)$.
\end{lem}
\begin{pf}
By Proposition \ref{pps4.1}, $D=\{u_{i,2}: i\equiv 1 \,({\rm
mod}\,4)\}\cup \{u_{i,1},u_{i,3}: i\equiv 3 \,({\rm mod}\,4)\}\cup \{u_{n,2}\}$ 
is a minimum dominating set and by the symmetry of $G_{n,3}$,
$D'=\{u_{i,2}: i\equiv 0 \,({\rm mod}\,4)\}\cup \{u_{i,1},u_{i,3}:
i\equiv 2 \,({\rm mod}\,4)\}\{u_{1,2}\}$ is also a minimum dominating set.
It is clear that if we
delete any vertical edge in $G_{n,3}$ or any horizontal edge
$u_{i,1}u_{i+1,1}$ and $u_{i,3}u_{i+1,3}$ where $i\equiv 0, 1\ or\
3\,({\rm mod}\,4)$ or any horizontal edge $u_{i,2}u_{i+1,2}$ where
$i\equiv 1,2\ or\ 3\,({\rm mod}\,4)$, $D$ or $D'$ is also a
domination set. Next, we consider the domination number of
$G_{n,3}-e$ where $e$ is an any other edge.

Let $e=u_{i,1}u_{i+1,1}$ or $e=u_{i,3}u_{i+1,3}$ where $i\equiv
2\,({\rm mod}\,4)$, or $e=u_{i,2}u_{i+1,2}$ where $i\equiv 0\,({\rm mod}\,4)$. 
Then $D''=\{u_{i,1},u_{i,3}:i\equiv 1\,({\rm mod}\,4)\}\cup
\{u_{i,2}:i\equiv 3\,({\rm mod}\,4)\}\cup \{u_{n,2}\}$ is 
a dominating set of $G-e$ with cardinality $n-\lfloor\frac{n-1}{4}\rfloor$.
By Lemma~\ref{lem2.1}, $|D''|=\gamma(G_{n,3})$.

From the above discussions, $\gamma(G_{n,3}-e)=\gamma(G_{n,3})$ for
any edge $e\in E(G_{n,3})$. Thus $b(G_{n,3} )\geqslant 2$.
\end{pf}

\begin{lem}\label{lem4.5}
$b(G_{n,3})\geqslant 2$ for $n\equiv 3\,({\rm mod}\,4)$.
\end{lem}

\begin{pf}
By Proposition \ref{pps4.1}, $D=\{u_{i,2}: i\equiv 1 \,({\rm
mod}\,4)\}\cup \{u_{i,1},u_{i,3}: i\equiv 3 \,({\rm mod}\,4)\}$ is a
minimum dominating set and by the symmetry of $G_{n,3}$,
$D'=\{u_{i,2}: i\equiv 3 \,({\rm mod}\,4)\}\cup \{u_{i,1},u_{i,3}:
i\equiv 1 \,({\rm mod}\,4)\}$ is also a minimum dominating set. It
is clear that if we delete any edge from $G_{n,3}$, $D$ or $D'$ is
also a dominating set. Thus $b(G_{n,3} )\geqslant 2$.
\end{pf}

Summing the above results, we have the following theorem,
immediately.

\begin{thm}\label{thm4.1} For $n\geqslant 3$,
$b(G_{n,3})=\left\{ \begin{array}{l}
 1,\ {\rm if}\ n\equiv 1\ {\rm or}\ 2 \,({\rm mod}\,4)  \\
2, \ {\rm if}\ n\equiv 0\ {\rm or}\ 3\,({\rm mod}\,4).
\end{array}
 \right. $
\end{thm}

\section{The bondage number of $G_{n,4}$}

In this section, let $A=\{1,2,3,5,6,9\}$.

\begin{lem}\label{lem5.1}
Let $D$ be a minimum dominating set of $G_{n,4}$. Then $1\leq
|Y_1\cap D|\leqslant 2$ and $1\leqslant |Y_n\cap D|\leqslant 2$
for $n\notin A$.
\end{lem}

\begin{pf}
By Lemma~\ref{lem2.1}, $|D|=n$. First, we prove that $1\leq |Y_1\cap
D|\leqslant 2$ and $1\leqslant |Y_n\cap D|\leqslant 2$. By the
symmetry of $G_{n,4}$, we only need to prove that $1\le |Y_1\cap
D|\leqslant 2$. By contradiction. Suppose $|Y_1\cap D|=0$ or
$|Y_1\cap D|\geqslant 3$.

If $|Y_1\cap D|=0$, then $|Y_2\cap D|=4$. By Lemma~\ref{lem2.2},
$|D\cap V(H_{n-2,4})|\geqslant \gamma(H_{n-3,4})$. By
Lemma~\ref{lem2.1}, $|D|\geqslant 4+\gamma(H_{n-3,4})\geq
4+n-3=n+1$, a contradiction with $|D|=n$. Thus $|Y_1\cap D|\geqslant
1$.

Assume now $|Y_1\cap D|\geqslant 3$. By Lemma~\ref{lem2.2}, $|D\cap
V(H_{n-1,4})|\geqslant \gamma(H_{n-2,4})$. By Lemma~\ref{lem2.1},
$|D|\geqslant 3+\gamma(H_{n-2,4})\geqslant 3+n-2=n+1$, a
contradiction with $|D|=n$. Thus $|Y_1\cap D|\leqslant 2$.
\end{pf}

\begin{lem}\label{lem5.2}
Let $D$ be a minimum dominating set of $G_{n,4}$. Then $|Y_1\cap
D|=1$, $|Y_n\cap D|=1$ for $n\in\{4,7,8,10,11\}$.
\end{lem}

\begin{pf}
By the symmetry of $G_{n,4}$ and by Lemma \ref{lem5.1},
we only need to prove $|Y_1\cap D|\neq 2$. 
Suppose, to the contrary, that that there exists a
minimum dominating set $D$ of $G_{n,4}$ such that $|Y_1\cap D|=2$.

If $n\ne 10$ then, by Lemma~\ref{lem2.2}, $|D\cap
V(H_{n-1,4})|\geqslant \gamma(H_{n-2,4})$. By Lemma \ref{lem2.1},
$|D|\geqslant 2+\gamma(H_{n-2,4})\ge 2+n-1=n+1$, a contradiction with
$|D|=n$.

Now assume $n=10$. Let $D'=D\setminus Y_1$. If $Y_2\cap D\neq
\emptyset$, then there exists a vertex $u_{2,j}$ such that $D'\cup
\{u_{2,j}\}$ is a dominating set of $H_{n-1,4}$. By Lemma
\ref{lem2.1}, $|D|=2+|D'| \geqslant 2+\gamma(G_{9,4})-1=11$, a
contradiction with $|D|=10$. Next, we assume that $Y_2\cap
D=\emptyset$ and then $|Y_3\cap D|\geqslant 2$. By
Lemma~\ref{lem2.2}, $|D\cap V(H_{n-3,4})|\geqslant
\gamma(H_{n-4,4})$. By Lemma \ref{lem2.1}, $|D|\geqslant
4+\gamma(H_{n-4,4})=4+7=11$, a contradiction with $|D|=10$. The
Lemma follows.
\end{pf}

\begin{lem}\label{lem5.3}
Let $D$ be a minimum dominating set of $G_{n,4}$. Then $|Y_1\cap
D|=1$, $|Y_n\cap D|=1$ for $n\notin A$.
\end{lem}

\begin{pf}
By the symmetry of $G_{n,4}$ and by Lemma \ref{lem5.1}, we only need
to prove $|Y_1\cap D|\neq 2$. By Lemma \ref{lem5.2}, the statement 
is true for $n\in \{4,7,8,10,11\}$. We proceed by induction on $n\geqslant
12$. 

Suppose that the assertion is true for any integer $k$ with
$10\leqslant k<n$. Suppose, to the contrary, that there exists a
minimum dominating set $D$ of $G_{n,4}$ such that $|Y_1\cap D|=2$. 
If $Y_2\cap D=\emptyset$, then $D'=D\setminus Y_1$ is a dominating set
of $H_{n-2,4}$ and $|Y_3\cap D'|\geqslant 2$. By the induction
hypothesis, $D'$ is not a minimum dominating set of $H_{n-2,4}$, and
hence $|D'|\geqslant \gamma(H_{n-2,4})+1\ge n-1$ by Lemma~\ref{lem2.1}.
Then $|D|=2+|D'|\geqslant n+1$, a contradiction with $|D|=n$.

If $Y_2\cap D\neq \emptyset$, there exists a vertex $u_{2,j}$ such
that $D''=(D\setminus Y_1)\cup \{u_{2,j}\}$ is a dominating set of
$H_{n-1,4}$ and $|Y_2\cap D''|\geqslant 2$. By the induction
hypothesis, $D''$ is not a minimum dominating set of $H_{n-1,4}$, and
hence $|D''|\geqslant \gamma(H_{n-1,4})+1\ge n$. Then $|D|\geqslant
2+|D''|-1\geqslant n+1$, a contradiction with $|D|=n$. The Lemma
follows.
\end{pf}

\begin{thm}\label{thm5.1}
$b(G_{5,4})=b(G_{9,4})=3$, $b(G_{6,4})=2$, and $b(G_{n,4})=1$ for $n\notin A$.
\end{thm}
\begin{pf}
It is easy to verify that $b(G_{5,4})=b(G_{9,4})=3$ and $b(G_{6,4})=2$.
Next, we prove $b(G_{n,4})=1$ for $n\notin A$. Then $n\ge 4$.
Let $D$ be a minimum dominating set of $G_{n,4}-u_{1,2}u_{1,3}$. By
Lemma~\ref{lem2.1}, we only need to show that $|D|\geqslant n+1$. We
prove the conclusion by considering three cases, respectively.

\begin{description}
\item[Case 1] $|Y_1\cap D|=0$.

Then $|Y_2\cap D|=4$. By Lemma~\ref{lem2.2}, $|D\cap
V(H_{n-2,4})|\geqslant \gamma(H_{n-3,4})$. Thus $|D|\geq
4+\gamma(H_{n-3,4})\ge n+1$.

\item[Case 2] $|Y_1\cap D|\geqslant 2$.

Then $D$ is a dominating set of $G_{n,4}$ with $|Y_1\cap D|\geqslant
2$. By Lemma \ref{lem5.3}, $D$ is not a minimum dominating set of
$G_{n,4}$, and hence $|D|\geqslant n+1$ by Lemma \ref{lem2.1}.

\item[Case 3] $|Y_1\cap D|=1$.

Without loss of generality, let $u_{1,j_0}\in D$ and $j_0\leqslant 1$. 
Then $u_{2,3},u_{2,4}\in D$ and hence $|Y_2\cap D|\geqslant 2$.

Let $D'=D\setminus\{u_{1,j_0}\}$. If $j_0=2$, or $|Y_2\cap D|\ge 3$,
or $j_0=1$ and $u_{31}\in D'$, then $D'$ is a dominating set of $H_{n-1,4}$
and let $D''=D'$.
Assume now $j_0=1$, $u_{31}\notin D$, and $Y_2\cap D=\{u_{2,3},u_{2,4}\}$.
If $u_{3,2}$ or $u_{3,3}$ or $u_{3,4}$ belongs to
$D$, then $D''=(D'\setminus \{u_{2,3}\})\cup \{u_{2,2}\}$ is a
dominating set of $H_{n-1,4}$ with $|Y_2\cap D''|\ge 2$. 

If $n\in \{4,7,10\}$, then $|D''|\ge \gamma(H_{n-1,4})=n$ by Lemma~\ref{lem2.1}.
If $n\notin \{4,7,10\}$,
then $D''$ is not a minimum dominating set of $H_{n-1,4}$
by Lemma~\ref{lem5.3}.
By Lemma~\ref{lem2.1}, $|D''|\ge \gamma(H_{n-1,4})+1=n$.
Thus $|D|\ge |D''|+1\ge n+1$.

In the following assume $j_0=1$, $u_{31}\notin D$, $Y_2\cap D=\{u_{2,3},u_{2,4}\}$,
and $x_{3,2},u_{3,3},u_{3,4}\notin D$. Then $u_{4,1},u_{4,2}$ should
be in $D$ to dominate $u_{31}$ and $u_{32}$
and $D'''=D\setminus \{u_{1,1},u_{2,3},u_{2,4}\}$ is a dominating set of
$H_{n-3,4}$ with $|Y_3\cap D'''|\ge 2$.

If $n\in \{4,8,12\}$, then $|D'''|\ge \gamma(H_{n-3,4})=n-2$ by Lemma~\ref{lem2.1}.
If $n\notin \{4,8,12\}$, then $D'''$ is not a minimum dominating set of
$H_{n-3,4}$ by Lemma~\ref{lem5.3}. Therefore
$|D'''|\ge \gamma(H_{n-3,4})+1=n-2$ by Lemma~\ref{lem2.1}.
Thus $|D|\ge 3+|D'''|\ge n+1$

\end{description}
The theorem follows.
\end{pf}

\section{Remarks}

Through determining the bondage number of $G_{n,m}$ for $2\leq m
\leqslant 4$, we find that if we delete the vertex $u_{1,1}$, the
domination number is invariable. If $m$ increases, the effect of
$u_{1,1}$ for the domination number will be smaller and smaller in
view of probability. Therefore we expect that $\gamma
(G_{n,m}-u_{1,1})=\gamma(G_{n,m})$ for $m\geqslant 5$ and we give
the following conjecture.

\begin{cjt}\label{cjt5.1}
$b(G_{n,m})\leqslant 2$ for $m\geqslant 5$.
\end{cjt}

In our method, determining the bondage number of a graph strongly
depends on the domination number of the graph. Even the domination
number of some graph, determining its bondage number is also very
difficult. For example, the domination number of $G_{n,m}$ for $m=5$
or $6$ has been determined~\cite{cch93,cch94}, we can not determined
its bondage number in our method since the cases that we need to
consider are much too. Thus, if we want to determine the bondage
number of $G_{n,m}$ for $m\geqslant 5$ or to solve the
Conjecture~\ref{cjt5.1}, we need a new method except for determining
the domination number of $G_{n,m}$ for $m\geqslant 7$. It is what we
further study work.

\end{document}